\documentclass[12pt]{amsart}

\usepackage{amssymb}
\usepackage{amsmath}
\evensidemargin\oddsidemargin

\begin{document}

\setcounter{page}{1}

\newtheorem{PROP}{Proposition}
\newtheorem{REMS}{Remark}
\newtheorem{LEM}{Lemma}
\newtheorem{THE}{Theorem\!\!}
\renewcommand{\theTHE}{}

\newcommand{\eqnsection}{
\renewcommand{\theequation}{\thesection.\arabic{equation}}
    \makeatletter
    \csname  @addtoreset\endcsname{equation}{section}
    \makeatother}
\eqnsection

\def\a{\alpha}
\def\b{\beta}
\def\B{{\bf B}} 
\def\CC{{\mathbb{C}}} 
\def\cia{c_{\a, \infty}}
\def\coa{c_{\a, 0}}
\def\cua{c_{\a, u}}
\def\cL{{\mathcal{L}}} 
\def\Ea{E_\a}
\def\EE{{\mathbb{E}}} 
\def\g{{\gamma}} 
\def\G{{\bf \Gamma}} 
\def\i{{\rm i}}
\def\K{{\bf K}}
\def\Ka{{\bf K}_\a}
\def\L{{\bf L}}
\def\lbd{\lambda}
\def\lcr{\left[}
\def\lpa{\left(}
\def\lva{\left|}
\def\pb{{\mathbb{P}}}
\def\rl{{\mathbb{R}}}
\def\prst{{\prec_{st}}}
\def\prcvx{{\prec_{cx}}}
\def\rpa{\right)}
\def\rcr{\right]}
\def\rva{\right|}
\def\W{{\bf W}}
\def\X{{\bf X}}
\def\U{{\bf U}}
\def\Un{{\bf 1}}
\def\Z{{\bf Z}}

\def\claw{\stackrel{d}{\longrightarrow}}
\def\elaw{\stackrel{d}{=}}
\def\qed{\hfill$\square$}
                  
\title[Unimodality of hitting times for stable processes]
      {Unimodality of hitting times for stable processes}

\author[Julien Letemplier]{Julien Letemplier}

\address{Laboratoire Paul Painlev\'e, Universit\'e Lille 1, Cit\'e Scientifique, F-59655 Villeneuve d'Ascq Cedex. {\em Email} : {\tt ju.letemplier@gmail.com}}

\author[Thomas Simon]{Thomas Simon}

\address{Laboratoire de physique th\'eorique et mod\`eles statistiques, Universit\'e  Paris Sud, B\^atiment 100, 15 rue Georges Cl\'emenceau, F-91405 Orsay Cedex. {\em Email} : {\tt simon@math.univ-lille1.fr}}



\begin{abstract} We show that the hitting times for points of real $\a-$stable L\'evy processes ($1<\a\le 2$) are unimodal random variables. The argument relies on strong unimodality and several recent multiplicative identities in law. In the symmetric case we use a factorization of Yano et al. \cite{Y3}, whereas in the completely asymmetric case we apply an identity of the second author \cite{TS1}. The method  extends to the general case thanks to a fractional moment evaluation due to Kuznetsov et al. \cite{KKMW}. 
\end{abstract}

\maketitle
 
\section{Introduction and statement of the result}

A real random variable $X$ is said to be unimodal if there exists $a\in\rl$ such that its distribution function $\pb[X\le x]$ is convex on $(-\infty, a)$ and concave on $(a, +\infty)$. When $X$ is absolutely continuous, this means that its density is non-decreasing on $(-\infty, a]$ and non-increasing on $[a, +\infty)$. The number $a$ is called a mode of $X$ and might not be unique. A random variable with a single mode is called strictly unimodal. The problem of unimodality has been intensively studied for infinitely divisible random variables and we refer to Chapter 10 in \cite{S} for details. This problem has also been settled in the framework of hitting times of processes and R\"osler - see Theorem 1.2 in \cite{R} - showed that hitting times for points of real-valued  diffusions are always unimodal. However, much less is known when the underlying process has jumps, for example when it is a L\'evy process. 

In this paper we consider a real strictly $\a-$stable process ($1<\a\le 2$), which is a L\'evy process $\{X_t, \; t\ge 0\}$ starting from zero and having characteristic exponent
$$\log[\EE[e^{\i \lbd X_1}]]\; =\; -(\i \lbd)^\a e^{-\i\pi\a\rho\, {\rm sgn}(\lbd)}, \quad \lbd\in\rl,$$
where $\rho\in[1-1/\a, 1/\a]$ is the positivity parameter of $\{X_t, \; t\ge 0\}$ that is $\rho =\pb[X_1 \ge 0].$ We refer to \cite{Z} and to Chapter 3 in \cite{S} for an account on stable laws and processes. In particular, comparing the parametrisations (B) and (C) in the introduction of \cite{Z} shows that the characteristic exponent of $\{X_t, \; t\ge 0\}$ takes the more familiar form
$$c\,\vert\lbd\vert^\a (1 - \i\theta\tan(\frac{\pi\a}{2})\,{\rm sgn}(\lbd))$$
with $\rho = 1/2 + (1/\pi\a) \tan^{-1}(\theta\tan(\pi\a/2))$ and $c = \cos(\pi\a(\rho -1/2)).$ The constant $c$ is a scaling parameter which could take any arbitrary positive value without changing our purposes below. We are interested in the hitting times for points of $\{X_t, \; t\ge 0\}\!\!:$ 
$$\tau_x\; =\; \inf\{t>0, \; X_t = x\}, \quad x\in\rl.$$
It is known  \cite{M} that $\tau_x$ is a proper random variable which is also absolutely continuous. Recall also - see e.g. Example 43.22 in \cite{S} - that points are polar for strictly $\a-$stable L\'evy process with $\a\le 1,$ so that $\tau_x = +\infty$ a.s. in this situation. In the following we will focus on the random variable $\tau =\tau_1.$ Again, this does not cause any loss of generality since by self-similarity one has 
$$\tau_x \,\elaw \,x^\a \tau_1\qquad\mbox{and}\qquad \tau_{-x}\,\elaw\, x^\a \tau_{-1}$$ 
for any $x \ge 0,$ and because the law of $\tau_{-1}$ can be deduced from that of $\tau_1$ in considering the dual process $\{-X_t, \; t\ge 0\}.$ We show the

\begin{THE} The random variable $\tau$ is unimodal.
\end{THE}

In the spectrally negative case $\rho = 1/\a,$ the result is plain because $\tau$ is then a positive stable random variable of order $1/\a,$ which is known to be unimodal - see e.g. Theorem 53.1 in \cite{S}. We will implicitly exclude this situation in the sequel and focus on the case with positive jumps. To proceed with this non-trivial situation we use several facts from the recent literature, in order to show that $\tau$ factorizes into the product of a certain unimodal random variable and a product of powers of Gamma random variables. The crucial property that the latter product preserves unimodality by independent multiplication \cite{CT} allows to conclude. Our argument follows that of \cite{TS2}, where a new proof of Yamazato's theorem for the unimodality of stable densities was established, but it is more involved. In passing we obtain a self-decomposability property for the Kanter random variable, which is interesting in itself and extends the main result of \cite{PSS}.

For the sake of clarity we divide the proof into three parts. We first consider the symmetric case $\rho =1/2$, appealing to a factorization of $\tau$ in terms of generalized Rayleigh and Beta random variables which was discovered in \cite{Y3}. Second, we deal with the spectrally  positive case $\rho = 1-1/\a$, with the help of a multiplicative identity in law for $\tau$ involving positive stable and shifted Cauchy random variables which was obtained in \cite{TS1}. In the third part, we observe the remarkable fact that this latter identity extends to the general case $\rho \in(1-1/\a, 1/\a)$, thanks to the evaluation of the Mellin transform of $\tau$ which was performed in \cite{KKMW}, and which can actually be obtained very easily - see section 2.4. This identity allows also to show that the density of $\tau$ is real-analytic and that $\tau$ is strictly unimodal - see the final remark (a).

\section{Proof of the Theorem}

\subsection{The symmetric case} Formula (5.12) and Lemma 2.17 in \cite{Y3} yield  together with the normalization (4.1) therein, which is the same as ours, the following independent factorization
\begin{equation}
\label{3Y}
\tau\; \elaw\;2^{-\a} \L^{-\frac{\a}{2}}\;\times\; \lpa\Z_{\frac{\a}{2}}^{(-\frac{1}{2})}\rpa^{-\frac{\a}{2}}\,\times\,\B_{1-\frac{1}{\a}, \frac{1}{\a}}^{-1}
\end{equation}
with the following notation, which will be used throughout the text:\\

\begin{itemize}

\item $\L$ is the unit exponential random variable.

\item $\B_{a,b}$ is the Beta random variable ($a,b >0$) with density 
$$\frac{\Gamma(a+b)}{\Gamma(a)\Gamma(b)} \, x^{a-1} (1-x)^{b-1} \Un_{(0,1)}(x)$$

\item  $\Z_c$ is the positive $c-$stable random variable ($0<c\le 1$), normalized such that
$$\EE\lcr e^{-\lbd \Z_c}\rcr \; =\; e^{-\lbd^c}, \qquad \lbd \ge 0.$$

\item For every $t\in\rl$ and every positive random variable $X$ such that $\EE[X^t]<+\infty,$ the random variable $X^{(t)}$ is the size-biased sampling of $X$ at order $t$, which is defined by
$$\EE\lcr f(X^{(t)})\rcr\; =\; \frac{\EE\lcr X^t f(X)\rcr}{\EE\lcr X^t\rcr}$$
for every $f : \rl^+\to\rl$ bounded continuous. 

\end{itemize}

\medskip

Notice that the above random variable $\Z_{\frac{\a}{2}}^{(-\frac{1}{2})}$ makes sense from the closed expression of the fractional moments of $\Z_c$:
\begin{equation}
\label{MomZ}
\EE[\Z_c^s]\; =\; \frac{\Gamma(1-s/c)}{\Gamma(1-s)}, \qquad s < c,
\end{equation}
which ensures $\EE[\Z_{\frac{\a}{2}}^{-\frac{1}{2}} ]< +\infty.$ Observe also that for every $t > 0$ one has $\L^{(t-1)}\elaw \G_t$ where $\G_t$ is the Gamma random variable with density
$$\frac{x^{t-1}e^{-x}}{\Gamma(t)} \, \Un_{(0,+\infty)}(x).$$
It is easy to see that
\begin{equation}
\label{Size}
X^{(t)}\times\,Y^{(t)}\; \elaw\; (X\times Y)^{(t)}\qquad\mbox{and}\qquad \lpa X^{(t)}\rpa^p\; \elaw\; \lpa X^p\rpa^{(\frac{t}{p})}
\end{equation}
for every $t,p\in\rl$ such that the involved random variables exist, and where the products in the first identity are supposed to be independent. In particular, one has
$$(\kappa X)^{(t)}\; \elaw\; \kappa X^{(t)}$$
for every positive constant $\kappa.$ Combined with (\ref{3Y}) the second identity in (\ref{Size}) entails
$$\tau\; \elaw\;2^{-\a} \L^{-\frac{\a}{2}}\;\times\; \lpa\Z_{\frac{\a}{2}}^{-\frac{\a}{2}}\rpa^{(\frac{1}{\a})}\,\times\,\B_{1-\frac{1}{\a}, \frac{1}{\a}}^{-1}.$$
On the other hand, Kanter's factorization - see Corollary 4.1 in \cite{K} - reads
\begin{equation}
\label{Kant}
\Z_{\frac{\a}{2}}^{-\frac{\a}{2}}\; \elaw\; \L^{1-\frac{\a}{2}}\,\times\; b_{\frac{\a}{2}}(\U)
\end{equation}
where $\U$ is uniform on $(0,1)$ and 
$$b_c (u)\; = \;\frac{\sin (\pi u)}{\sin^c (\pi c u)\sin^{1-c} (\pi(1-c) u)}$$
for all $u,c\in (0,1).$ Since $b_c$ is a decreasing function from $\kappa_c = c^{-c}(1-c)^{c-1}$ to $0$ - see the proof of Theorem 4.1 in \cite{K} for this fact, let us finally notice that the support of the random variable
$$\K_{c}\; =\; \kappa_{c}^{-1} b_{c}(\U)$$ 
is $[0,1].$ Putting everything together shows that
$$\tau\; \elaw\;2^{-\a}  \kappa_{\frac{\a}{2}}^{}\, \L^{-\frac{\a}{2}}\times\; \lpa\L^{1-\frac{\a}{2}}\rpa^{(\frac{1}{\a})}\,\times\, \K_{\frac{\a}{2}}^{(\frac{1}{\a})}\,\times\,\B_{1-\frac{1}{\a}, \frac{1}{\a}}^{-1}.$$
By Theorem 3.7. and Corollary 3.14. in \cite{CT}, the random variable 
$$\X_\a\; =\; 2^{-\a}  \kappa_{\frac{\a}{2}}^{}\, \L^{-\frac{\a}{2}}\times\; \lpa\L^{1-\frac{\a}{2}}\rpa^{(\frac{1}{\a})}\; \elaw\; 2^{-\a}  \kappa_{\frac{\a}{2}}^{}\, \L^{-\frac{\a}{2}}\times\; \lpa\L^{(\frac{1}{\a}-\frac{1}{2})}\rpa^{1-\frac{\a}{2}}\; \elaw\; 2^{-\a}  \kappa_{\frac{\a}{2}}^{}\, \L^{-\frac{\a}{2}}\times\; \G_{\frac{1}{\a}+\frac{1}{2}}^{1-\frac{\a}{2}}$$
is multiplicatively strongly unimodal, that is its independent product with any unimodal random variable remains unimodal. Indeed, a straightforward computation shows that the random variables $\log (\L)$ and $\log ( \G_{\frac{1}{\a}+\frac{1}{2}})$ have a log-concave density, and the same is true for $\log(\X_\a)$ by Pr\'ekopa's theorem.
All in all, we are reduced to show the

\begin{PROP}
\label{Ksym}
With the above notation, the random variable $\K_{\frac{\a}{2}}^{(\frac{1}{\a})}\times\B_{1-\frac{1}{\a}, \frac{1}{\a}}^{-1}$ is unimodal.
\end{PROP}

\proof A computation shows that the density of $\B_{1-\frac{1}{\a}, \frac{1}{\a}}^{-1}$  decreases on $(1,+\infty).$ Hence, there exists  $F_\a : (0,1)\mapsto (1,+\infty)$  increasing and convex such that
$$\B_{1-\frac{1}{\a}, \frac{1}{\a}}^{-1}\;\elaw\; F_\a(\U).$$
On the other hand, up to normalization the density of $\K_{\frac{\a}{2}}^{(\frac{1}{\a})}$ writes
$$g_\a(x)\; =\;x^{\frac{1}{\a}} f_\a(x)$$
on $(0,1),$ where $f_\a$ is the density of $\K_{\frac{\a}{2}}.$ It follows from Lemma 2.1 in \cite{TS3} that $f_\a$ increases on $(0,1),$ so that the density of $g_\a$ also increases on $(0,1)$ and that there exists  $G_\a : (0,1)\mapsto (0,1)$  increasing and concave such that
$$\K_{\frac{\a}{2}}^{(\frac{1}{\a})}\;\elaw\; G_\a(\U).$$
We can now conclude by the lemma in \cite{TS2}.

\endproof

\begin{REMS} {\em The lemma in \cite{TS2} shows that the mode of $\K_{\frac{\a}{2}}^{(\frac{1}{\a})}\times\B_{1-\frac{1}{\a}, \frac{1}{\a}}^{-1}$ is actually 1. However, this does not give any information on the mode of $\tau.$}
\end{REMS}

\subsection{The spectrally positive case} This situation corresponds to the value $\rho = 1-1/\a$ of the positivity parameter. The characteristic exponent of $\{X_t, \; t\ge 0\}$ can be extended to the negative half-plane, taking the simple form
$$\log[\EE[e^{-\lbd X_1}]]\; =\; \lbd^\a, \qquad \lbd \ge 0.$$
With this normalization, we will use the following independent factorization which was obtained in \cite{TS1}:
$$\tau\; \elaw\; \U_\a\;\times\; \Z_{\frac{1}{\a}}$$
where $\U_\a$ is a random variable with density
$$f_{\U_\a}(t)\; =\;\frac{-(\sin \pi\a) t^{1/\a}}{\pi(t^2 - 2t \cos \pi\a  +1)}\cdot$$
It is easy to see that $\U_\a$ is multiplicatively strongly unimodal for $\a\le 3/2$ and this was used in Proposition 8 of \cite{TS1} to deduce the unimodality of $\tau$ in this situation. To deal with the general case $\a\in (1,2)$ we proceed via a different method. First, it is well-known and easy to see - solve e.g. Exercise 4.21 (3) in \cite{CY} - that the independent quotient
$$\lpa\frac{\Z_{\a-1}}{\Z_{\a-1}}\rpa^{\a-1}$$
has the density
$$\frac{-\sin (\pi \a)}{\pi (t^2 -2 t\cos \pi \a + 1)}$$
over $\rl^+,$ whence we deduce
\begin{eqnarray*}
\tau & \elaw & \lpa\frac{\Z_{\a-1}^{\a-1}}{\Z_{\a-1}^{\a-1}}\rpa^{(\frac{1}{\a})}\times\; \Z_{\frac{1}{\a}}\\
& \elaw & \kappa_{\frac{1}{\a}}^{-\a}\lpa\frac{\L^{2-\a}}{\L^{2-\a}}\rpa^{(\frac{1}{\a})}\times\; \L^{1-\a}\times\;\K_{\a-1}^{(\frac{1}{\a})}\times\; (\K_{\a-1}^{-1})^{(\frac{1}{\a})}\times\K_{\frac{1}{\a}}^{-\a}
%\\& \elaw & \kappa_{\frac{1}{\a}}\lpa\frac{\G_{\frac{2}{\a}-1}}{\G_{\frac{2}{\a}-1}}\rpa^{2-\a}\times\; \L^{1-\a}\times\;\K_{\a-1}^{(\frac{1}{\a})}\times\; (\K_{\a-1}^{-1})^{(\frac{1}{\a})}\times\K_{\frac{1}{\a}}^{-\a}
\end{eqnarray*}
with the above notation. Similarly as above, the first product with the three exponential random variables is multiplicatively strongly unimodal, whereas the random variable $\K_{\a-1}^{(\frac{1}{\a})}$ has an increasing density on $(0,1).$ Hence, reasoning as in Proposition \ref{Ksym} it is enough to show that the random variable $(\K_{\a-1}^{-1})^{(\frac{1}{\a})}\times\K_{\frac{1}{\a}}^{-\a}$ has a decreasing density on $(1, +\infty).$ We show the more general

\begin{PROP}
\label{Kasym}
With the above notation, the random variable
$$(\K_{\b}^{-r})^{(t)}\,\times\,\K_{\g}^{-s}$$
has a decreasing density on $(1, +\infty)$ for every $r, s >0$ and $\b, \g, t$ in $(0,1).$
\end{PROP}

The proof of the proposition uses the notion of self-decomposability - see Chapter 3 in \cite{S} for an account. Recall that a positive random variable $X$ is self-decomposable if its Laplace transform reads
$$\EE[e^{-\lbd X}]\; =\; \exp - \lcr a_X^{}\lbd \; +\; \int_0^\infty (1- e^{-\lbd x}) \frac{\varphi_{X}^{}(x)}{x} dx\rcr, \quad \lbd\ge 0,$$
for some $a_X^{}\ge 0$ which is called the drift coefficient of $X,$ and some non-increasing function $\varphi_{X}^{} : (0,+\infty)\to \rl^+$ which will be henceforth referred to as the spectral function of $X.$ Introduce the following random variable
$$\W_\b\; =\; -\log(\K_\b)$$
and notice that its support is $\rl^+,$ thanks to our normalization for $\K_\b.$ A key-observation is the following

\newpage

\begin{LEM}
\label{Main}
The random variable $\W_\b$ is self-decomposable, without drift and with a spectral function taking the value $1/2$ at $0\!+\!.$
\end{LEM}

\proof Combining (\ref{MomZ}), (\ref{Kant}) and the classical formula for the Gamma function
$$\Gamma(1-u)\; =\; \exp\lcr \gamma u \, +\, \int_0^\infty (e^{ux}-1-u x)\frac{dx}{x(e^x-1)}\rcr, \quad u<1,$$
(where $\gamma$ is Euler's constant) yields the following expression for the Laplace transform of $\W_\b$ - see (3.5) in \cite{PSS}:
$$\EE[e^{-\lbd \W_\b}]\; =\; \EE[\K_\b^{\lbd}]\; =\; \exp - \lcr \int_0^\infty (1- e^{-\lbd x}) \frac{\varphi_{\b}^{}(x)}{x} dx\rcr, \quad \lbd\ge 0,$$
with
$$\varphi_\b(x)\; =\; \frac{e^{-x}}{1-e^{-x}}\; -\; \frac{e^{-x/\b}}{1-e^{-x/\b}}\; -\; \frac{e^{-x/(1-\b)}}{1-e^{-x/(1-\b)}}, \quad x > 0.$$
It was shown in Lemma 3 of \cite{PSS} that the function $\varphi_{\b}$ is non-negative and an asymptotic expansion at order 2 yields $\varphi_{\b}(0+) = 1/2.$ 

We finally show that $\varphi_{\b}$ is non-increasing on $(0,+\infty).$ Following the proof and the notation of Lemma 3 in \cite{PSS} this amounts to the fact that the function $x \mapsto x\psi_\b(x)$ therein is non-decreasing on $(0,1),$ which is a clear consequence of the following claim
\begin{equation}
\label{Clay}
t\; \mapsto\; \log(1-e^t) - \log (1-e^{rt})\quad \mbox{is convex on $\rl^-$ for every $r\in (0,1).$}
\end{equation}
Let us show the claim. Differentiating twice, we see that we are reduced to prove that
$$\frac{r^2x^{r-1}(1-x)^2}{(1-x^r)^2}\; \ge \; 1, \qquad 0<x,r< 1.$$
The limit of the quantity on the left-hand side is 1 when $x\to1-,$ whereas its logarithmic derivative equals
$$\frac{(r-1)(1-x^{r+1}) + (r+1)(x^r -x)}{x(1-x)(1-x^r)}\cdot$$
The latter fraction is negative for all $0<x,r<1$ because its numerator is concave as a function of $x\in(0,1)$ which vanishes together with its derivative at $x=1\!-\!.$ This shows (\ref{Clay}) and finishes the proof of the lemma.

\endproof

\noindent
{\bf Proof of Proposition \ref{Kasym}.} Set $f_{\b, r}$ resp. $f_{\g,s}$ for the density of $r\W_\b$ resp. $s\W_\g.$ By multiplicative convolution, the density of $(\K_{\b}^{-r})^{(t)}\,\times\,\K_{\g}^{-s}$ writes
$$\int_1^x (xy^{-1})^t f^{}_{\K_{\b}^{-r}} (xy^{-1}) f^{}_{\K_{\g}^{-s}} (y) \frac{dy}{y}$$
on $(1, +\infty),$ up to some normalization constant. This transforms into
$$x^{t-1}\int_1^x f_{\b, r}(\log (x) -\log(y)) f_{\g,s}(\log(y))\frac{dy}{y^{t+1}}\; =\; x^{t-1}\int_0^{\log(x)} f_{\b, r}(\log (x) - u) e^{-tu} f_{\g,s}(u)\, du$$
and since $t\in(0,1)$ it is enough to prove that the function
\begin{equation}
\label{Yam}
v\;\mapsto\; \int_0^v f_{\b, r}(v-u) e^{-tu} f_{\g,s}(u)\, du
\end{equation}
is non-increasing on $(0,+\infty).$ Lemma \ref{Main} and a change of variable show that $f_{\b, r}(u)$ resp. $e^{-tu} f_{\g,s}(u)$ is up to normalization the density of a positive self-decomposable random variable without drift and with spectral function $\varphi_{\b}(xr^{-1})$ resp. $e^{-tx} \varphi_{\g}(xs^{-1}),$ with the notation of the proof of Lemma \ref{Main}. By additive convolution this entails that the function in (\ref{Yam}) is the constant multiple of the density of a positive self-decomposable random variable without drift and with spectral function 
$$\varphi_{\b}(xr^{-1})\; +\; e^{-tx} \varphi_{\g}(xs^{-1}).$$
By Lemma  \ref{Main} this latter function takes the value 1 at $0+,$ and we can conclude by Theorem 53.4 (ii) in \cite{S}.

\qed

\begin{REMS} {\em By Theorem 4 in \cite{JS} we know that $\W_\b$ has also a completely monotone density, in other words - see Theorem 51.12 in \cite{S} - that its spectral function writes
$$\varphi_\b (x)\; =\; x\int_0^\infty\!\! e^{-tx} \theta_\a (t)dt, \quad x\ge 0,$$
for some function $\theta_\a (t)$ valued in $[0,1]$ and such that
$t^{-1}\theta_\a (t)$ is integrable at $0+\!.$ This entails that $\K_{\b}^{-r}$ has a completely monotone density as well for every $r>0$ - see Corollary 3 in \cite{JS}. However, this latter property does not seem true in general for $(\K_{\b}^{-r})^{(t)}\,\times\,\K_{\g}^{-s}.$}
\end{REMS}

\subsection{The general case} We now suppose $\rho\in (1-1/\a, 1/\a),$ which means that our stable L\'evy process has jumps of both signs. The symmetric case was dealt with previously but it can also be handled with the present argument. Theorem 3.10 in \cite{KKMW} computes the fractional moments of $\tau$ in closed form:
\begin{equation}
\label{closed}
\EE[\tau^s]\; =\; \frac{\sin(\frac{\pi}{\a})\sin(\pi\rho\a(s+\frac{1}{\a})) }{\sin(\pi\rho)\sin(\pi(s+\frac{1}{\a}))}\; \times\;\frac{\Gamma(1-\a s)}{\Gamma(1-s)}, \quad -1-\frac{1}{\a} < s < 1-\frac{1}{\a}
\end{equation}
(the initial normalization of \cite{KKMW} is the same as ours - see the introduction therein - but beware that with their notation our $\tau$ has the law of $T_0$ under ${\rm P}_{\!-1}$). On the other hand, it is easy to see from (\ref{MomZ}) and the complement formula for the Gamma function that 
$$\EE\lcr \lpa\frac{\Z_{\rho\a}^{\rho \a}}{\Z_{\rho\a}^{\rho \a}}\rpa^{s}\rcr\; =\; \frac{\sin(\pi\rho\a s)}{\rho\a\sin(\pi s)},\quad -1< s< 1.$$
Hence, a fractional moment identification entails
\begin{equation}
\label{Rk}
\tau\; \elaw\; \lpa\frac{\Z_{\rho\a}^{\rho \a}}{\Z_{\rho\a}^{\rho \a}}\rpa^{(\frac{1}{\a})}\times\; \Z_{\frac{1}{\a}}\cdot
\end{equation}
Making the same manipulations as in the spectrally positive case, we obtain
\begin{equation}
\label{Final}
\tau \;\elaw\;\kappa_{\frac{1}{\a}}^{-\a}\lpa\frac{\L^{1-\rho\a}}{\L^{1-\rho\a}}\rpa^{(\frac{1}{\a})}\times\; \L^{1-\a}\times\;\K_{\rho\a}^{(\frac{1}{\a})}\times\; (\K_{\rho\a}^{-1})^{(\frac{1}{\a})}\times\K_{\frac{1}{\a}}^{-\a}
\end{equation}
and we can conclude because Proposition \ref{Kasym} applies here as well. 

\qed

\subsection{A short proof of (\ref{closed})} In this paragraph we give an independent proof of the fractional moment evaluation (\ref{closed}), which is short and standard. This method was suggested to us by L.~Chaumont and P.~Patie and we would like to thank them for reminding us this classical argument. Setting $f_{X_t}$ for the density of $X_t,$ by Theorem 43.3 in \cite{S} one has
$$\EE[e^{-q\tau}]\; =\; \frac{u^q(1)}{u^q(0)}$$
where by self-similarity
$$u^q(1)\; =\; \int_0^\infty e^{-qt} f_{X_t}(1)\, dt \; =\; \int_0^\infty e^{-qt} t^{-1/\a}f_{X_1}(t^{-1/\a})\, dt,$$
and
$$u^q(0)\; =\; \int_0^\infty e^{-qt} f_{X_t}(0)\, dt \; =\; \frac{1}{\varphi(q)}$$
with $\varphi$ the Laplace exponent of the inverse local time at zero of $\{X_t, \, t\ge 0\}.$ It is well-known and easy to see by self-similarity - see e.g. Theorem 2 in \cite{St} - that
$$\varphi(q) \; =\; \kappa\, q^{\frac{\a-1}{\a}},$$
where $\kappa > 0$ is a normalizing constant to be determined later. Making a change of variable, we deduce
$$\EE[e^{-q\tau}]\; =\; {\tilde \kappa}\, q^{\frac{\a-1}{\a}} \EE[Z_1 e^{-q Z_1}]$$
with ${\tilde \kappa} =\a\rho\kappa$ and $Z_1 = (X_1\vert X_1\ge 0)^{-\a}.$ For every $s\in (0,1)$ one has
\begin{eqnarray*}
\EE[\tau^{-s}] & = & \frac{1}{\Gamma(s)} \int_0^\infty \EE[e^{-q\tau}] q^{s-1} dq\\
& = & \frac{ {\tilde \kappa}}{\Gamma(s)} \EE\lcr Z_1 \int_0^\infty e^{-q Z_1} q^{\frac{\a-1}{\a} +s-1} dq\rcr\\
& = & {\tilde \kappa}\, \frac{\Gamma(1-1/\a+s)}{\Gamma(s)}\, \EE\!\lcr Z_1^{\frac{1}{\a} -s} \rcr \; =\;  {\tilde \kappa}\, \frac{\Gamma(1-1/\a+s)}{\Gamma(s)}\, \EE\!\lcr (X_1\vert X_1\ge 0)^{\a s -1} \rcr.\\
\end{eqnarray*}
On the other hand, by the formula (2.6.20) in \cite{Z} one has
$$\EE\!\lcr (X_1\vert X_1\ge 0)^{\a s -1} \rcr\; =\; \frac{\Gamma(\a s ) \Gamma(1+1/\a-s)}{\Gamma(1-\rho +\rho\a s ) \Gamma(1+\rho -\rho\a s)}$$
and putting everything together entails
\begin{eqnarray*}
\EE[\tau^{-s}] & = & \frac{{\tilde \kappa} \Gamma(\a s ) \Gamma(1-1/\a+s) \Gamma(1+1/\a-s)}{\Gamma(s)\Gamma(1-\rho +\rho\a s ) \Gamma(1+\rho -\rho\a s)}
\; = \;\frac{\sin(\frac{\pi}{\a})\sin(\pi\rho\a(-s+\frac{1}{\a})) \Gamma(1+\a s)}{\sin(\pi\rho)\sin(\pi(-s+\frac{1}{\a}))\Gamma(1+s)}
\end{eqnarray*}
for every $s\in (0,1),$ where we used standard properties of the Gamma function and the identification of the constant comes from $\EE[\tau^{0}] =1.$ This completes the proof of (\ref{closed}) for $s\in (-1,0),$ and hence for $s\in (-1-1/\a, 1- 1/\a)$ by analytic continuation. 

\qed

\begin{REMS}{\em The above computation shows that the normalizing constant for the inverse local time at zero reads
$$\kappa\; =\; \frac{\a\sin(\frac{\pi}{\a})}{\sin(\pi\rho)}\cdot$$ 
One can check that this constant is the same as the one computed by Fourier inversion in \cite{St} p. 636.} 
\end{REMS}

\subsection{Final remarks} (a) The identity in law (\ref{Rk}) shows that the density of $\tau$ is the multiplicative convolution of two densities which are real-analytic on $(0,+\infty).$ Indeed, it is well-known - see e.g. \cite{Z} Theorem 2.4.1 - that the density of $\Z_{\frac{1}{\a}}$ is real-analytic on $(0,+\infty),$ whereas the density of the first factor in (\ref{Rk}) reads
$$\frac{\sin (\pi \rho\a)\sin(\frac{\pi}{\a})t^{\frac{1}{\a}}}{\pi \sin(\pi\rho)(t^2 +2 t\cos(\pi \rho\a) + 1)}\cdot$$
Hence, the density of $\tau$ is itself real-analytic on $(0,+\infty)$ and a combination of our main result and the principle of isolated zeroes entails that $\tau$ is strictly unimodal. Besides, its mode is positive since we know from Theorem 3.15 (iii) in \cite{KKMW} in the spectrally two-sided case, and from Proposition 2 in \cite{TS1} in the spectrally positive case, that the density of $\tau$ always vanishes at $0+$ (with an infinite first derivative). The strict unimodality of $\tau$ can also be obtained in analyzing more sharply the factors in (\ref{Final}) and using Step 6 p. 212 in \cite{CT}.\\

(b) The identity in law (\ref{Rk}) can be extended in order to encompass the whole set of admissible parameters $\{\a\in (1, 2],\,1-1/\a \le \rho\le 1/\a\}$ of strictly $\a-$stable L\'evy processes that hit points in finite time a.s. Using the Legendre-Gauss multiplication formula and a fractional moment identification, one can also deduce from (\ref{Rk}) with $\rho =1/2$ the formula (5.12) of \cite{Y3}. It is possible to derive a factorization of $\tau$ with the same inverse Beta factor for $\a \le n/(n-1)$ and $\rho =1/n,$ but this kind of identity in law does not seem to be true in general.\\

(c) When $\rho\a \ge 1/2,$ the formula (\ref{Rk}) shows that the law of $\tau$ is closely related to that of the positive branch of a real stable random variable with scaling parameter $1/\a$ and positivity parameter $\rho\a.$ Indeed, Bochner's subordination - see e.g. Chapter 6 in \cite{S} or Section 3.2 in \cite{Z} for details - shows that the latter random variable decomposes into the independent product  
$$\lpa\frac{\Z_{\rho\a}^{\rho \a}}{\Z_{\rho\a}^{\rho \a}}\rpa\times\, \Z_{\frac{1}{\a}}\cdot$$

(d) The identity (\ref{Rk}) is attractive in its simplicity. Compare with the distribution of first passage times of stable L\'evy processes, whose fractional moments can be computed in certain situations - see Theorem 3 in \cite{Kz} and the references therein for other recent results in the same vein - but with complicated formul\ae\, apparently not leading to tractable multiplicative identities in law. In the framework of hitting times, it is natural to ask whether (\ref{Rk}) could not help to investigate further distributional properties of $\tau,$ in the spirit of \cite{Y}. This will be the matter of further research.

\medskip
  
\noindent
{\bf Acknowledgement.}  Ce travail a b\'en\'efici\'e d'une aide de l'Agence Nationale de la Recherche portant la r\'ef\'erence ANR-09-BLAN-0084-01.

\end{document}